\documentclass{article}
\usepackage{amsmath,amssymb}
\usepackage[latin1]{inputenc}
\usepackage[dvips]{graphics}
\input{diagrams.tex}

\newcommand{\Supp}{\operatorname{Supp}}
\newcommand{\Hom}{\operatorname{Hom}}
\newcommand{\Ext}{\operatorname{Ext}}
\newcommand{\End}{\operatorname{End}}
\newcommand{\add}{\operatorname{add}}
\newcommand{\lra}{\leftrightarrow}
\newcommand{\mra}{minimal right $\textup{add} \overline{T}$-approximation
  in $\C$}
\newcommand{\mla}{minimal left $\textup{add} \overline{T}$-approximation
  in $\C$}

\newcommand{\QT}{\mathcal{Q}_T}
\newcommand{\qt}{\mathcal{Q}_{\overline{T}}}
\newcommand{\QC}{\mathcal{Q}_C}
\newcommand{\qc}{\mathcal{Q}_{\overline{C}}}
\newcommand{\IC}{{I}_C}
\newcommand{\IT}{{I}_T}
\newcommand{\ic}{{I}_{\overline{C}}}
\newcommand{\Qalt}{\mathcal{Q}_{\textup{alt}}}
\newcommand{\F}{F}

\newcommand{\C}{\mathcal{C}}

\newcommand{\ssi}{\Leftrightarrow}

\newcommand{\ot}{\leftarrow}
\newcommand{\za}{\alpha}
\newcommand{\zb}{\beta}
\newcommand{\ze}{\varepsilon}
\newcommand{\zg}{\gamma}

\newtheorem{theorem}{Theorem}[section]
\newtheorem{proposition}[theorem]{Proposition}
\newtheorem{conjecture}[theorem]{Conjecture}
\newtheorem{corollary}[theorem]{Corollary}
\newtheorem{lemma}[theorem]{Lemma}

\newtheorem{definition}{Definition}

\newtheorem{remark}[theorem]{Remark}
\newtheorem{example}[theorem]{Example}

\newenvironment{proof}{\begin{trivlist}\item{\bf{Proof.}}}
  {\hfill\rule{2mm}{2mm}\end{trivlist}}

\title{Quivers with relations and cluster tilted algebras}
\author{Philippe Caldero, Fr\'ed\'eric Chapoton, Ralf Schiffler}
\date{}
\begin{document}
\maketitle
\begin{abstract}
  Cluster algebras were introduced by S. Fomin and A. Zelevinsky in
  connection with dual canonical bases.
  To a cluster algebra of simply laced Dynkin type one can associate the
  cluster category. Any cluster of the cluster algebra corresponds to a
  tilting object in the cluster category. The cluster
  tilted algebra is the algebra of endomorphisms of that tilting object.
  Viewing the cluster tilted algebra as a path algebra of a quiver with
  relations, we prove in this paper that the quiver of the cluster
  tilted algebra  is equal to the cluster diagram.
  We study also the  relations. As an application of
  these results, we answer several conjectures on the connection
  between cluster algebras and quiver representations.
\end{abstract}

\setcounter{section}{-1}


\begin{section}{Introduction}
Cluster algebras were introduced in the work of S. Fomin
and A. Zelevinsky, \cite{cluster1,cluster2,Ysystems}.
This theory appeared in the context of dual canonical basis and more
particularly in the study of the Berenstein-Zelevinsky conjecture.
Cluster algebras are now connected with many topics: double Bruhat
cells, Poisson varieties, total positivity, Teichmüller spaces. The
main results on cluster algebras are on the one hand the
classification of finite cluster algebras by root systems and on the
other hand the realization of algebras of regular functions on double
Bruhat cells in terms of cluster algebras.

Recently, many new results have been established relating cluster
algebras of simply laced finite type to quiver representations.
It has been shown in \cite{CCS} (type $A$) and \cite{BMRRT} (types
$A,D,E$) that the set of cluster variables is
in bijection with the set of indecomposable objects in the so called
cluster category $\C$, which is the quotient category $D/\tau^{-1}[1]$
of the bounded derived category $D$ of quiver representations  by the inverse Auslander-Reiten
translate $\tau^{-1}$ composed with the shift $[1]$.

For type $A$  the
authors associated in \cite{CCS} a quiver with relations  to each
cluster 
in such a
way that the indecomposable representations of that quiver with
relations are in bijection with all cluster variables outside the
cluster. A result of this approach was the description of the
denominator of the Laurent polynomial expansion of any cluster variable in
the variables of any cluster. In this paper, we generalize this result
to the types $D$ and $E$ (Theorem \ref{denominator}).

Buan, Marsh, Reineke, Reiten and Todorov \cite{BMRRT}
used tilting theory to relate the cluster algebra to the cluster
category; each cluster corresponds to a tilting object in $\C$.
For several concepts in the theory of cluster algebras,
they obtained nice module theoretic interpretations, e.g. exchange
pairs, compatibility degree. They called the endomorphism algebra of
their tilting object \textit{cluster tilted algebra} and conjectured
that this algebra is isomorphic to the path algebra of our quiver with
relations \cite[Conj. 9.2]{BMRRT}. In this paper we prove this conjecture in type $A$
(Theorem \ref{A}) and parts of it in types $D$ and $E$,
(Theorem \ref{quivers} and Proposition \ref{relations}).
We also prove another of their conjectures \cite[Conj. 9.3]{BMRRT} on the module theoretic
calculation of the exchange relations in the cluster algebra,
(Theorem \ref{exchange}).
Buan, Marsh and Reiten also announced results on these conjectures.
 In \cite{BMR}, Buan, Marsh and Reiten studied further the cluster tilted
algebra and gave a precise description of its module category.

The paper is organized as follows. In section \ref{cta}  we recall briefly some
facts about cluster tilted algebras. Lemma \ref{BB'} is a new result,
but it follows almost immediately from \cite{BMRRT}. In section \ref{ca} we
list some concepts of cluster algebras that we will need later.
Section \ref{qar} contains the crucial results. We prove there that the quiver
of the cluster tilted algebra  is equal to the cluster diagram and that
the relations defined in \cite{CCS} are also satisfied in that algebra.
In section \ref{appl}, we prove the conjectures mentioned above. They
follow easily from the results in section \ref{qar}. Finally, in the
Appendix we include some general results on embeddings of
cluster diagrams in the plane.
\end{section}

\begin{section}{Cluster tilted algebras} \label{cta}
Let $k$ be an algebraically closed field and
$\Qalt$  an alternating quiver of simply-laced Dynkin type,
$D$ the bounded derived category of
finitely generated modules with shift functor $[1]$ and $\C=D/\tau^{-1}[1]$
the cluster category of \cite{BMRRT}.
Here $\tau$ denotes the
Auslander-Reiten translate.
By a result of Keller \cite{Keller}, $\C$ is a triangulated category.
Let $P_1,\ldots,P_n$ be the indecomposable projective $k\Qalt$-modules
and $I_1,\ldots,I_n$  the injective ones.
According to \cite{BMRRT},
there is a natural fundamental domain of indecomposable objects for
$\C=D/\tau^{-1}[1]$ in $D$ consisting of $\{\textup{indecomposable
$k\Qalt$-modules}\}\cup\{P_i[1]\mid i=1,\ldots,n\}$.
We will think of the indecomposable objects of $\C$ as their
representatives in this fundamental domain. Thus an indecomposable object
$M$ in
$\C$ is either a $k\Qalt$-module or $M=P_l[1]$ for some $l$.
An indecomposable object $M$ in $\C$ is called $l$-free if it satisfies
$l\notin \Supp M$ if $M$ is a module and $l\ne i$ if $M=P_i[1]$.
Throughout this paper we will use the notation
$[M,N]_A=\dim\Hom_A(M,N)$ and $  [M,N]_A^1=\dim\Ext_A(M,N)$,
for $A=k\Qalt, \C ,D$.

Let $C$ be a cluster of the cluster algebra of the same type as $\Qalt$,
and let $T=\oplus_{i=1}^n T_i$ be the corresponding tilting object  of the
cluster category $\C$ \cite{BMRRT}.
The following lemma is proved in
\cite[Lemma 8.2]{BMRRT}.
\begin{lemma}\label{hom}
  $[T_i,T_j]_\C\le 1$.
\end{lemma}
Let $(\mathcal{Q}_T,I_T)$ be a quiver  with
relations such that its path algebra $k\mathcal{Q}_T /\langle I_T\rangle$ is
isomorphic to  the cluster tilted algebra $\End_\C(T)^{op}$ of \cite{BMR}.
Hence the
vertices of $\mathcal{Q}_T$ 'are' the indecomposable direct summands
$T_1,\ldots,T_n$ of $T$ and the lemma implies that
there is an arrow $T_j \to T_i$ precisely if
$[T_i,T_j]_\C=1$  and  no non-zero morphism $f\in \Hom_\C(T_i,T_j)$
factors through one of the $T_k, k\ne i,j$.

Following \cite{BMRRT},
let $\overline{T}=\oplus_{i=1}^{n-1}\, T_i$ be an almost complete basic tilting
object in $\C$ and let $M,M'$ be the two complements of $\overline{T}$. Then
$T=\overline{T}\oplus M$ and $T'=\overline{T}\oplus M'$ are tilting
objects, and there are triangles
  $M'\ \to \ B \ \to \ M\ \to \ M'[1]$
  and $M\ \to \ B' \ \to \ M'\ \to \ M[1]$ in $\C$,
  where $B \to M $ is a \mra~ and $M \to B'$ is a \mla.
  Recall that $f:B\to M$ (resp. $f':M\to B'$)
  being a minimal right (resp. left)
  $\add\overline{T}$-approximation in $\C$ means
  \begin{enumerate}
    \item $B$ (resp. $B'$) is an object of $\add\overline{T}$.
    \item The induced map  $\Hom_\C(X,B)\to \Hom_\C(X,M)$
    (resp. $\Hom_\C(B',X)\to \Hom_\C(M,X)$) is surjective for all objects
    $X$ of $\add\overline{T}$.
    \item (Minimality) For every map $g:B\to B$ (resp. $g':B'\to B'$) such that $fg=f$
    (resp. $g'f=f$), the map $g$ (resp. $g'$) is an
    isomorphism.
  \end{enumerate}
                             We have the following
\begin{lemma}\label{BB'}
$$B=\oplus_{i\in I} T_i,$$ where $I=\{i\mid M\to T_i \textup{ in }
\mathcal{Q}_{T}\}.$\\
$$B'=\oplus_{i\in I'} T_i,$$ where $I'=\{i\mid T_i\to M \textup{ in }
\mathcal{Q}_{T}\}.$
\end{lemma}

\begin{proof}
$B=\oplus_{i\in J}a_iT_i$, $a_i\ge 1$,
for some subset $J\subset \{1,2,\ldots,n-1\}$.
We will show first that all $a_i$ are equal to $1$.
Suppose $a_{l}>1$, write $B=X\oplus T_{l}\oplus \ldots \oplus T_{l}$
where $X$ has no direct summand isomorphic to $T_{l}$. Let $h$ be a
generator of
$\Hom_\C(T_l,M)$. Since $[T_l,M]_\C\le 1$, we can write
the minimal right \textup{add}$T$-approximation $f:B\to M$ in matrix form
as $f=[\hat f\
b_1h\ \ldots b_{a_l}h]$, for some scalars $b_1,\ldots,b_{a_l}$
and with $\hat f$  the restriction of $f$ to $X$.
If this matrix is $[\hat f\  0 \ldots 0]$ then take
$g_0=\left[\begin{array}{cc}\textup{Id}_X& 0\\ 0&0\end{array}\right]$ and
get $fg_0=f$. By minimality of $f$ we get that $g_0:B\to B$ is an
isomorphism, contradiction.
Otherwise  $h$ as well as one of the $b_i$ are non-zero. Say $b_2\ne 0$.
Put
$$g=\left[\begin{array}{cccccc}
\textup{Id}_X&0&0&\ldots&0\\
0&0&0&\ldots&0\\
0&b_1/b_2&1&\ldots&0\\
\vdots&\vdots&\vdots&\ddots&\vdots \\
0&0&0&\ldots&1
\end{array}\right]
$$
Then $fg=f$ and by minimality of $f$ we get that  $g$ is an
isomorphism, contradiction.

So $B=\oplus_{i\in J}T_i$.
To show that $I\subset J$, suppose that there exists $i_0$ such that
$M\to T_{i_0} $ in $\mathcal{Q}_{T}$ and $i_0 \notin J$; thus $T_{i_0}$ is
not a direct summand of $B$.
Since $M\to T_{i_0} $ in $\mathcal{Q}_{T}$, there is
a non-zero element $g$ in $\Hom_\C(T_{i_0},M)$.
By property 2 above, there exists $h:T_{i_0}\to B$
such that $g=fh$, and this implies that there is
no arrow $M\to T_{i_0}  $ in $\mathcal{Q}_{T'}$, contradiction.

To show that $J\subset I$, suppose that there is $i_0\in J$ such that
$M\to T_{i_0}$ is not in $\mathcal{Q}_{T}$.
Suppose first that $[T_{i_0},M]_\mathcal{C}=0$.
Hence the restriction of $f$ to $T_{i_0}$ is zero.
Write $B=\oplus_{i\in J\setminus i_0} T_i
\oplus T_{i_0}$ and define $g:B\to B$ in matrix block form to be
$$\left[\begin{array}{cc} 1&0\\0&0    \end{array}\right].$$
Then $fg=f$ and since $f:B\to M$ is minimal this
implies that $g$ is an isomorphism, contradiction.
Suppose now that $[T_{i_0},M]_\mathcal{C}=1$ and let $f_{i_0}:T_{i_0}\to M$ be
the restriction of $f$.
Then there exists an arrow $M\to T_l$ in $\mathcal{Q}_{T}$ such that
$[T_{i_0},T_l]_\C=[T_{l},M]_\C=1$, because
 $M\to T_{i_0}$ is not in $\mathcal{Q}_{T}$. Since we have already shown that
$I\subset J$, it is clear that $T_l$ is a direct summand of $B$.
Let $f_l:T_l \to M$ be  the restriction of $f$. If $f_l=0$, we get a contradiction as above.
Thus $f_l\ne 0$ and there exists $h:T_{i_0}\to T_l$ such that
$f_{i_0}=f_l h$, since the dimensions of all corresponding Hom-spaces is
1.
Write
$B=T_{i_0}\oplus T_l \oplus\left(\oplus_{i\in J\setminus
\{i_0,l\}} T_i \right)$
and define $g:B\to B$ in matrix block form to be
$$\left[\begin{array}{ccc} 0&0&0\\h&1&0\\0&0&1    \end{array}\right].$$
Note that $f=[f_{i_0}\,f_l\,f']$. Thus $fg=f$ and by minimality of $f$ we have
that $g$ is an isomorphism, contradiction.
The proof for $B'$ is similar and left to the reader.
\end{proof}
\end{section}


\begin{section}{Cluster algebras}\label{ca}
For the proof of Theorem \ref{quivers} we will need some concepts of
\cite{cluster2}. For convenience we recall them here briefly but our
exposition is only for simply laced finite types, i.e. $A,D,E$.

Let $I_+$ be the set of sinks of $\Qalt$ and $I_-$ the set of sources.
Define the sign function $\ze$ on vertices of $\Qalt$ by
$\ze(i)= +1$ if $i\in I_+$ and
$\ze(i)= -1$ if $i\in I_-$.
Let $\textup{Q}$ be the root lattice and $\Phi_{\ge -1}$ the set of almost
positive  roots. Denote the simple (positive) roots by
$\za_1,\ldots,\za_n$ and the corresponding simple reflections by
$s_1,\ldots,s_n$. Let $\tau_+,\tau_-$ be the involutions on $\Phi_{\ge
-1}$ given by
$$
\tau_\ze(\za)=\left\{\begin{array}{ll}
\za &\textup{if }\za=-\za_i ,\ i\in I_{-\ze}\\
(\prod_{i\in I_\ze}s_i)\za   &\textup{otherwise}
\end{array}\right.
\quad\quad \ze\in\{+,-\}
$$
Let $\langle \tau_+,\tau_-\rangle$  be the group generated by $\tau_+$
and $\tau_-$.
Note that the composition $\tau_-\circ\tau_+$ is the Coxeter
transformation on positive roots.

There is a bijection $\za\mapsto x_\za$ between the set of almost
positive roots and the set of cluster variables.
Two almost positive roots $\zb,\zb'$ are called exchangeable if there are two
clusters $C,C'$ such that $C'=C\setminus\{x_\zb\} \cup \{x_{\zb'}\} $.
\begin{proposition}\textup{\cite[Prop.3.3]{cluster2}}
Given  any $\zg\in \textup{Q}$, then there exists a cluster $C$ such
that $\zg$ can be written as
$\zg=\sum_{x_\za\in C} a_{\za\zg}x_\za $ with $a_{\za\zg}\ge 0$. The
almost positive roots $\za$ such that $a_{\za\zg}\ne 0$ in this
expansion are called \textup{cluster components} of $\zg$ with respect to the
cluster $C$.
\end{proposition}
\begin{proposition} \label{uplus} \textup{\cite[Prop.3.6]{cluster2}}
 If $\zb,\zb'\in \Phi_{\ge -1}$ are exchangeable then the set
 $$\{\sigma^{-1}(\sigma(\zb)+\sigma(\zb'))\mid \sigma \in
 \langle\tau_+,\tau_-\rangle\}$$
consists of two elements of $\textup{Q}$, one of
which is $\beta + \beta'$, and the other will be denoted by
$\beta \uplus \beta'$.
In the special case where $\zb'$ is the negative simple
  root $-\za_l$ we have
  $$\zb\uplus(-\za_l)=\zb-\za_l-\sum_{l\frac{\ \ }{\ \ } j \in \Qalt}
\za_j.$$
\end{proposition}

\begin{lemma}\label{sign}\textup{\cite[Lemma 4.1]{cluster2}}
There exists a sign function $(\beta, \beta') \mapsto
\varepsilon(\beta, \beta')\in \{\pm 1\}$ on pairs of exchangeable roots,
uniquely determined by the following properties:
$$\begin{array}{l}
\varepsilon(- \alpha_j,\beta') = - \varepsilon (j)\,;\\
 \varepsilon(\tau \beta,\tau \beta')= - \varepsilon(\beta,\beta')
\text{\ for $\tau\in\{\tau_+,\tau_-\}$
and $\beta,\beta' \notin
\{-\alpha_j \mid \tau(-\alpha_j) =-\alpha_j\}$.}
\end{array}$$
Moreover, this function is skew-symmetric:
$$\varepsilon (\beta',\beta) = - \varepsilon(\beta, \beta').$$
\end{lemma}

\end{section}


\begin{section}{Quivers and relations}\label{qar}
Thinking of the cluster tilted algebra as a path algebra of a quiver
with relations, we will prove in this section that the quiver in
question is the
cluster diagram. Moreover we will show that the relations defined in
\cite{CCS} for the cluster diagram are also satisfied in the cluster
tilted algebra.

 Let $(\QC,\IC)$ be the quiver with relations
associated to the cluster $C$ in \cite{CCS}.
Recall that
$\QC$ is the cluster diagram of the cluster $C$ as defined in \cite{cluster2}
and that the set of relations $I_C$ can be expressed as follows using
the notion of shortest paths.
By definition, a \textit{shortest path} in the quiver $\QC$ is an
oriented path (with no repeated arrow) contained in an induced
subgraph of $\QC$ which is a cycle.
For any arrow $i\to j$ in $\QC$, let $P_{ji}$
be the set of shortest paths from $j$ to $i$ in $\QC$.
We
will show in the Appendix, that for any arrow $i\to j$ the set
 $P_{ji}$ has at most 2 elements.
Define
$$p(j,i)=\left\{\begin{array}{ll}
p&\textup{if $P_{ji}=\{p\}$}\\
p_1-p_2&\textup{if $P_{ji}=\{p_1,p_2\}$.}
\end{array}\right. $$
Then
$$\IC=\bigcup_{i\to j} \big\{p(j,i)\big\}.$$
 Let      $\langle I_C\rangle$ be the ideal generated by $I_C$.

It has been conjectured in \cite[Conj. 9.2]{BMRRT} that $k\QC/\langle I_C\rangle$
is isomorphic
to the cluster tilted algebra $\End_\C(T)^{\textup{op}}$. We will show
that $\QC=\QT$ and
$\langle I_C\rangle\subset \langle I_T\rangle$.
In type $A$, we can then deduce the conjecture using
the fact that the number of indecomposable modules over both algebras is
equal.
\begin{theorem}\label{quivers} Let $C$ be any cluster of a cluster
algebra of type $A,D$ or $E$ and let $\QC$ be its cluster diagram. Let
$T$ be a corresponding tilting object in the cluster category and
$(\QT,I_T)$ the quiver with relations of the cluster tilted algebra
$\End_\C(T)^{\textup{op}}$. Then
$$\QT=\QC.$$
\end{theorem}

\begin{proof}
The vertices of $\QC$ are almost positive roots and will be denoted by
greek letters.
It has been shown in \cite[sect.3]{cluster2} that there is an arrow
$  \za\to\zb \textup{ in } \QC$ if and only if
  \begin{equation}\label{eq4}
  \left\{\begin{array}{ll}
  \textup{Either $\ze(\zb,\zb')=-1$ and $\za$ is a cluster component of
  $\zb \uplus \zb'$}\\
  \textup{or $\ze(\zb,\zb')=+1$ and $\za$ is a cluster component of
  $\zb + \zb'$,}
  \end{array}\right\}
  \end{equation}
  where $\zb' $ is the unique almost positive root such that $C\setminus
  \{\zb\}\cup \{\zb'\}$ is a cluster.

According to \cite{BMRRT}, to each almost positive root $\za$
corresponds an indecomposable object $M_\za$ in $\C$.
 Let $T=\overline{T}\oplus M_\zb$ and let
  $M_{\zb'}$ be the other complement of the almost complete basic tilting
  object $\overline{T}$. 
  The indecomposable object $M_{\zb'}$  in $\C$ corresponds   to  $\zb'$.
  We may suppose without loss of generality that $M_{\zb'}$ is the first
  shift of the $l$-th indecomposable projective module, $M_{\zb'}=P_l[1]$,
  for some $l$. That means that $\zb'=-\za_l$.
  Thus
  $$\ze(\zb,\zb')=\ze(\zb,-\za_l)=-\ze(-\za_l,\zb)=\left\{\begin{array}{ll}
  1 &\textup{if $l$ is a sink}\\
  -1&\textup{if $l$ is a source}\end{array}\right.
$$
where the last two identities follow from Lemma \ref{sign} and our
choice of the sign function $\ze$.
Let us suppose first that $M_\zb$ is different from $\tau M_{\zb'} $ and
$\tau^{-1} M_{\zb'}$, that is $M_\zb\ne I_l,P_l$.
Then by a result of $\cite{BMRRT}$
we have the following two triangles in $\C$
$$\begin{array}{ccccccc}
M_{\zb'}&\to& B &\to&M_\zb &\to&M_{\zb'}[1]\\
M_{\zb}&\to& B' &\to&M_{\zb'} &\to&M_{\zb}[1]\\
\end{array}$$
and $B=\oplus_{\zg\in I} M_\zg$ and
$B'=\oplus_{\zg\in I'} M_\zg$ with $I=\{\zg\mid M_\zb\to M_\zg \textup{ in
} \QT\}$
and $I'=\{\zg\mid M_\zb\ot M_\zg \textup{ in } \QT\}$, by Lemma \ref{BB'}.
Now using $M_{\zb'}=P_l[1]$ and $M_{\zb'}[1]=I_l$ these two triangles give
\begin{equation}\label{eq1}
\sum_{\zg\in I}\zg=\left\{\begin{array}{ll}
\zb-\dim I_l &\textup{if $l$ is a sink in $\Qalt$}\\
\zb-\za_l &\textup{if $l$ is a source in
$\Qalt$}\end{array}\right.\end{equation}
\begin{equation}\label{eq2}
\sum_{\zg\in I'}\zg=\left\{\begin{array}{ll}
\zb-\za_l &\textup{if $l$ is a sink in $\Qalt$}\\
\zb-\dim P_l &\textup{if $l$ is a source in
$\Qalt$}\end{array}\right.\end{equation}

{Indeed, let $B=B_0\oplus B_1$ with $B_0$ a $k\Qalt$-module and
$B_1=\oplus_{j\in J} P_j[1]$. Note that $B_1$ is zero if $l $ is a source in
$\Qalt$ by Lemma \ref{BB'},
and if $l$ is a sink then $J$ is a subset of the set of
neighbours of $l$ in $\Qalt$. In particular, all elements of $J$ are sources
in $\Qalt$ and thus the indecomposable injective  $I_j$ is a simple
module for $j\in J$.
Note also that $M_\zb$ is a $k\Qalt$-module since $[P_l[1],M_\zb]^1_\C\ne
0.$
With this notation, the first triangle gives the following
triangle in $D$:
$$I_l[-1] \to \bigoplus_{j\in J}I_j[-1]    \bigoplus B_0  \to M_\zb \to I_l.$$
We apply the functor $\Hom_D(P_i,-)$ to this triangle and get the following
exact sequence
$$\begin{array}{ccccccccccc}
0&\to &{\Hom_D(P_i,B_0)}&\to&\Hom_D(P_i,M_\zb)\\
&\to&\Hom_D(P_i,I_l)&\to&\Hom_D(P_i,\oplus_{j\in J}I_j) &\to&0
\end{array}$$
whence $\underline{\dim}(B_0)_i=\underline{\dim}(M_\zb)_i
-\underline{\dim}(I_l)_i +\sum_{j\in J}\delta_{ij}$
since  $I_j$ is a simple module for all $j\in J$.
Thus $\sum_{\zg\in I}\zg=\zb-\dim(I_l)$.
This implies equation (\ref{eq1}).
Now let $B'=B'_0\oplus B'_1$ with $B'_0$ a $k\Qalt$-module and
$B'_1=\oplus_{j\in J'} P_j[1]$. Note that $B'_1$ is zero if $l $ is a sink in
$\Qalt$, and if $l$ is a source then $J'$ is a subset of the set of
neighbours of $l$ in $\Qalt$. In particular, all elements of $J'$ are sinks
in $\Qalt$ and thus $P_j$ is a simple module for $j\in J'$.
The second triangle gives the following
triangle in $D$:
$$M_\zb \to \bigoplus_{j\in J'}P_j[1]    \bigoplus B'_0  \to P_l[1]\to
M_\zb[1] .$$
We apply the functor $\Hom_D(P_i,-)$ to this triangle and get the following
exact sequence
$$\begin{array}{ccccccccccc}
0&\to &\Hom_D(P_i,\bigoplus_{j\in J'}P_j) &\to&
{\Hom_D(P_i,P_l)}                    \\
&\to&\Hom_D(P_i,M_\zb)&\to&\Hom_D(P_i,B_0') &\to&0
\end{array}$$
whence $\underline{\dim}(B'_0)_i=\underline{\dim}(M_\zb)_i
-\underline{\dim}(P_l)_i +\sum_{j\in J'}\delta_{ij}$
since  $P_j$ is a simple module for all $j\in J'$.
Thus $\sum_{\zg\in I'}\zg=\zb-\dim(P_l)$.
This implies equation (\ref{eq2}).
}

By Proposition \ref{uplus},
$$\zb\uplus(-\za_l)=\zb-\za_l-\sum_{l\frac{\ \ }{\ \ } j \in \Qalt}
\za_j$$  and   $\za_l+\sum_{l\frac{\ \ }{\ \ } j \in \Qalt}\za_j$ is $\dim
P_l$ if $l$ is a source and $\dim I_l$ if $l$ is a sink.
Thus
\begin{equation}\label{eq6}
\sum_{\zg\in I}\zg=\left\{\begin{array}{ll}
\zb\uplus \zb' &\textup{if $l$ is a sink in $\Qalt$}\\
\zb+\zb' &\textup{if $l$ is a source in
$\Qalt$}\end{array}\right.\end{equation}
\begin{equation}\label{eq7}
\sum_{\zg\in I'}\zg=\left\{\begin{array}{ll}
\zb+\zb' &\textup{if $l$ is a sink in $\Qalt$}\\
\zb\uplus \zb' &\textup{if $l$ is a source in
$\Qalt$.}\end{array}\right.\end{equation}
Hence
if $l$ is a sink we have $\ze(\zb,\zb')=1$ and then
$$\begin{array}{rcl}\za\to\zb\textup{ in } \QC &\stackrel{(\ref{eq4})}{\Longleftrightarrow}& \za \textup{ is a cluster
component of $\zb + \zb'$}\\
&\stackrel{(\ref{eq7})}{\Longleftrightarrow}& M_\za \textup{ is a direct summand
 of $B'$}\\
&\stackrel{\textup{Lemma} \ \ref{BB'}}{\Longleftrightarrow}& M_\za \to M_\zb \textup{ in } \QT
\end{array}$$
and
$$\begin{array}{rcl}\za\ot\zb \textup{ in } \QC&\stackrel{(\ref{eq4})}{\Longleftrightarrow}& \za \textup{ is a cluster
component of $\zb\uplus \zb'$}\\
&\stackrel{(\ref{eq6})}{\Longleftrightarrow}& M_\za \textup{ is a direct summand
 of $B$}\\
&\stackrel{\textup{Lemma} \ \ref{BB'}}{\Longleftrightarrow}& M_\za \ot M_\zb  \textup{ in } \QT
\end{array}$$
and
if $l$ is a source we have $\ze(\zb,\zb')=-1$ and then
$$\begin{array}{rcl}\za\to\zb \textup{ in } \QC &\stackrel{(\ref{eq4})}{\Longleftrightarrow}& \za \textup{ is a cluster
component of $\zb \uplus \zb'$}\\
&\stackrel{(\ref{eq7})}{\Longleftrightarrow}& M_\za \textup{ is a direct summand
 of $B'$}\\
&\stackrel{\textup{Lemma}\  \ref{BB'}}{\Longleftrightarrow}& M_\za \to M_\zb \textup{ in } \QT
\end{array}$$
and
$$\begin{array}{rcl}\za\ot\zb \textup{ in } \QC &\stackrel{(\ref{eq4})}{\Longleftrightarrow}& \za \textup{ is a cluster
component of $\zb+ \zb'$}\\
&\stackrel{(\ref{eq6})}{\Longleftrightarrow}& M_\za \textup{ is a direct summand
 of $B$}\\
&\stackrel{\textup{Lemma} \ \ref{BB'}}{\Longleftrightarrow}& M_\za \ot M_\zb  \textup{ in } \QT
\end{array}$$

We still need to consider $M_\zb\in\{P_l,I_l\}$. These two cases are
similar and we will only treat the case $M_\zb=P_l$. Thus
$M_{\zb'}=\tau M_\zb$ in $\C$. Then there is
no minimal left $\textup{add} \overline{T}$-approximation $f':B'\to M_{\zb'}$
  in $\C$  and
we only have one triangle $M_{\zb'}\to B\to
M_{\zb}\to M_{\zb'}[1]$.
Furthermore, there is no arrow $M_\zg \to M_\zb $ in $\QT$ because otherwise
$1=[M_\zb,M_\zg]_\C=[M_\zg,\tau M_\zb]^1_\C=[M_\zg,M_{\zb'}]^1_\C$ which
contradicts the fact that $\overline{T}\oplus M_{\zb'}$ is a tilting
object.
Thus $M_\zb$  is  a sink in $\QT$.
Note that $\zb+\zb'=0$ if $l$ is a sink and $\zb\uplus \zb'=0$ if $l$ is a
source in $\Qalt$. Therefore there is no arrow $\za\to \zb$ in $\QC$ by
(\ref{eq4}).
On the other hand, Lemma \ref{BB'} still gives $B=\oplus_{\zg\in I} M_\zg$
and equations (\ref{eq1}) and (\ref{eq6}) as well as  the proof of the
equivalence $\za\ot\zb \ssi M_\za\ot M_\zb$ still hold as before.
This proves the theorem.
\end{proof}

Now we want to study the relations $I_T$. First we need to investigate
shortest paths. Let us write $\F$ for the composition $\tau^{-1}[1]$.
Given an indecomposable object $T$ in our fundamental
domain of $\C$, we say that
an indecomposable object $\widetilde{T}$ in $D$ is over $T$ if it lies in
the $\F$-orbit of $T$.
\begin{lemma}
Let $T \to T'$ be a non-zero morphism between indecomposable objects
in $\C$ and let $\widetilde{T}$ be an indecomposable object in $D$ over $T$.
   Then there exists a unique
indecomposable object $\widetilde{T}'$ in $D$ over $T'$ such that there is
a non-zero morphism $\widetilde{T} \to \widetilde{T}'$.
\end{lemma}
\begin{proof}
The existence of $\widetilde {T}'$ is clear since the $\C$-morphism $T\to
T'$ is non-zero.
Uniqueness follows easily from the well known fact that
for any two indecomposable objects $M,N$ in $D$ we
  have
  $$\left[M,N\right]_D \ne 0 \quad \Rightarrow \quad [M,N[a]]_D=0
  \textup{ for all $a\ne 0$}.$$
\end{proof}
Let $p:
T_1\stackrel{p_{12}}{\to} T_2\stackrel{p_{23}}{\to}\ldots
\stackrel{p_{(k-1)k}}{\to} T_k$
be a path in $\QT$ from $T_1$ to $T_k$,
that is a composition of arrows $T_i\stackrel{p_{i(i+1)}}{\to}T_{i+1}$.
Denote by $p^{\C}$ the corresponding element of $\End_\C(T)^{\textup{op}}$
under the isomorphism $k\QT/\langle I_T
\rangle\cong\End_\C(T)^{\textup{op}}$.
Then $p^{\C}\in\Hom_\C(T_k,T_1)$ is the composition of morphisms
$p^\C=p_{12}^\C\circ p_{23}^\C\circ \ldots \circ p_{(k-1)k}^\C$ with
each $ p_{i(i+1)}^\C\in\Hom_\C(T_i,T_{i+1})$ non-zero.
Let us construct a \textit{lift} $p^D$
 of $p$ as follows.
 Consider first the case where $p$ is a single arrow $T_1 \to T_2$.
 Then by the lemma, given an indecomposable
 object $\widetilde {T}_2$ in $D$ over $T_2$ there exists a unique indecomposable
 object $\widetilde {T}_1$ over $T_1$   such that there is
a non-zero morphism $\widetilde{T}_2 \to \widetilde{T}_1$. Any such non-zero
morphism is called a \textit{lift of $T_1\to T_2$ starting at
$\widetilde{T}_2$}.
Note that this lift is
unique up to multiplication by a  scalar.
Now let $p$ be  any path.
 We choose an indecomposable object $\widetilde{T}_k$
 over $T_k$.
 Using the lemma on each morphism $p_{i(i+1)}^\C$, there is
 a unique family of indecomposable objects $(\widetilde {T}_i)_
 {i=k-1,\ldots,1}$,
 with  $\widetilde {T}_i$  over $T_i$, and  such that
  $  \Hom_D(\widetilde {T}_{i+1},\widetilde{ T}_i)\ne 0$.
For each arrow $T_i\to T_{i+1}$ let $p_{i(i+1)}^D$ be a lift of
$p_{i(i+1)}$ starting at $\widetilde{T}_{i+1}$.
Then the
    composition of morphisms
    $p^D=p^D_{12}\circ p^D_{23}\circ \ldots p^D_{(k-1)k} $
    is called a \textit{lift of $p$ starting at $T_k$}.
Note that
  $p^D\in\Hom_D(\widetilde{T}_k,\widetilde{T}_1)$ is unique up to multiplication
  by scalar. Note also that $p^D$ may be zero although each
 $p_{i(i+1)}^D$ is non-zero.

\begin{definition} If  $p$ is a closed path
 in the situation above, that is $T_1=T_k$, then
    $p^D\in\Hom_D(\widetilde{T}_1,\F^{a} \widetilde{T}_1)$ and $a$ is called the
    \textup{winding number} of the path $p$.
\end{definition}
\begin{proposition}\label{shortest} Let $p:T_1\to T_2\to \ldots \to
T_{k}=T_1$ be a closed path and let $p_{21}$ be the
    subpath $T_2\to T_3\to \ldots \to T_k=T_1$.
  \begin{itemize}
    \item[(1)] The winding number $a$ is zero if and only if $p$ is a
    constant path.
    \item[(2)] Suppose that $p$ is not constant.
    Then $p_{21} $ is a shortest path if the winding number
    $a$ is equal to $1$.
    \item[(3)] If $p_{21}$ is not a shortest path then $p_{21}$ is zero
    in $(\QT,I_T)$.
  \end{itemize}
\end{proposition}
\begin{proof}
(1) is obvious. To show (2), suppose $a=1$.
If $p_{21}$ is not a shortest path, then the subquiver of $\QT$ induced by
the vertices on $p_{21}$ is not a cycle. Now since cycles in cluster diagrams
are always oriented, there exist 2 vertices $T_i,T_j$ on the path
$p_{21}$ such that $p_{21}=p_{2i}\,p_{ij}\,p_{j1}$, where $p_{kl}$ is a path from
$T_k$ to $T_l$, and such that there is an arrow
$T_i\stackrel{\kappa}{\to} T_j$ in $\QT$ and the subpath $p_{ij}$ of $p$
is not an arrow.
Then the path $\kappa \,p_{j1}\, p_{1i} $ is a non-constant closed path in
$\QT$, thus its winding number is at least $1$.
Since $p_{ij}$ is not an arrow, it follows that the winding number of
the path $p_{ij} \,p_{j1}\, p_{1i} $ is at least 2, contradiction.
This proves (2).

Suppose now that (3) is not true. That is, $p_{21}$ is a non-zero
non-shortest path. Suppose without loss of generality that $T_2$ is an
indecomposable projective $k\Qalt$-module.
Since $p_{12}$ is an arrow in $\QT$, its lift $p^D_{12}$ is
non-zero and thus $T_1$ is an indecomposable
$k\Qalt$-module too.
On the other hand, $p_{21}^D\in  \Hom_D(T_1,\tau^{-a} T_2[a])$
  is non-zero, so $a=1$ and $p_{21} $ is a shortest path by (2). This
  proves (3).
\end{proof}
\begin{conjecture} Suppose the situation of Proposition
\ref{shortest}(2).
 Then  the winding number $a$ is equal to $1$
  if $p_{21} $ is a shortest path.
    \end{conjecture}

\begin{proposition}\label{relations}
Let $C$ be any cluster of a cluster
algebra of type $A,D$ or $E$ and let $(\QC,I_C)$ be the associated
quiver with relations. Let
$T$ be a corresponding tilting object in the cluster category and
$(\QT,I_T)$ the quiver with relations of the cluster tilted algebra
$\End_\C(T)^{\textup{op}}$. Then
  $$\langle\IC\rangle\subset \langle\IT\rangle.$$
\end{proposition}

\begin{proof}  Let $T_j\to T_i$ be an arrow
  in $\QT$ and $\mathcal{P}_{ij}=\{p_1,p_2,\ldots,p_m\}$ be the set of shortest paths from
  $T_i$  to $T_j$ in $\QT$. We have to show that
  $$\begin{array}{cc}
  p_1=0 &\textup{if }m=1 \\
  p_1=p_2 &\textup{if }m=2 \\
  \end{array}$$
    The proof is by induction on the rank $n$.
    The smallest case is $n=3$. In this case $\QT$ is either a Dynkin
    quiver of type $A_3$ and then $m=0$ or $\QT$ is the cyclic quiver of
    rank $3$.
    In the latter case we may suppose without loss of generality that $T_j$
    is the $l$-th indecomposable projective $k\Qalt$-module $P_l$.
    Then $l$ is a leaf of $\Qalt$, $T_k=I_l$ the $l$-th
    indecomposable injective module and $T_i=P_{l'}[1]$, where $l'$ is the
    other leaf of $\Qalt$. We illustrate this situation in the
    Auslander-Reiten quiver of $\C$ in Figure \ref{A3}. Obviously $[T_j,T_i]_\C=0$, whence $p_1:T_i\to
    T_k\to T_j$ is a zero path in $\QT$.

    \begin{figure}
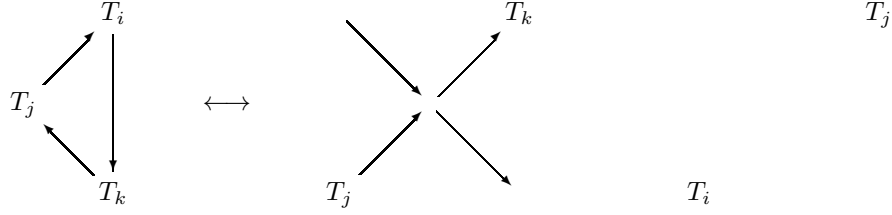

    \begin{diagram}[size=1.7em]
        &    &T_i&       &&&&&T_k&&&&&&&&T_j\\
       &\ruTo&&&&\rdTo&&\ruTo&&\rdTo&&\ruTo&&\rdTo&&\ruTo\\
      T_j &&\dTo&\hspace*{1cm}\longleftrightarrow\hspace*{1cm}&& &         \\
         &\luTo&&&&\ruTo&&\rdTo&&\ruTo&&\rdTo&&\ruTo&&\rdTo\\
         &   &T_k&&T_j&&&&&&&&T_i
      \end{diagram}
      \caption{Cyclic quiver of rank $3$ and corresponding Auslander-Reiten quiver with tilting object}\label{A3}
    \end{figure}

From now on let $n>3$.
  We will show the case $m=1$ first.
  For convenience, let us relabel the vertices of $\QT$ such that
 $p=p_1:T_1\to T_2\to \ldots \to
  T_k$.
  Suppose there exists a $T_{i_0}$ such that the path $p$ does not pass
  through   $T_{i_0}$, i.e. $k<n$.
  We may suppose without loss of generality that
  $T_{i_0}=P_l[1]$ for some $l$.
  Let   $\overline{C}=C\setminus\{-\za_l\}$.
  By a result of \cite{cluster2}, $\overline{C}$ is a cluster of a cluster
  algebra of rank $n-1$ and its quiver with relation $(\qc,\ic)$
  is the full subquiver $\qc$ of $\QC$ with vertices
  $\QC^0\setminus\{-\za_l\}$ and $\ic$ its usual set of relations. Note
  that $\qc$ may be disconnected.
  Let $\overline{\C}$ be the cluster category of the quiver $\Qalt
  \setminus{l}$ and denote by $\overline{T_h}$ ($h\ne i_0$)
  the restriction of the $l$-free object $T_h$ to $\overline{\C}$.
  Let $\overline{T}=\oplus_{h\ne i_0}\overline{ T_h}$. Then
  $\overline{T}$ is the tilting object in
  $\overline{\C}$  that corresponds to the cluster $\overline{C}$.
  We have already shown in Theorem \ref{quivers}
  that $\qc=\qt$ and by induction we conclude that the path $\overline{p}$,
  which is the "restriction" of the path $p$ to $\qt$,
  is zero in $(\qt,I_{\overline{T}})$. We want to show that $p$ is zero in
  $(\QT,\IT)$.
  Suppose the contrary.
That is
\begin{equation}\label{eq10}
[\overline{T_k},\overline{T_1}]_{\overline{\C}}=0 \quad
\textup{and} \quad [T_k,T_1]_\C=1 .
\end{equation}
Let us show first that $[T_k,T_1]_D=0$.
We will proceed using a case by case analysis. $T_1$ (resp. $T_k$) may
be either an indecomposable $k\Qalt$-module or the first shift of an
indecomposable  projective $k\Qalt$-module different from
$T_{i_0}=P_l[1]$, hence there are 4 different
cases to consider.
\begin{enumerate}
\item  $T_1$ and $T_k$ are $k\Qalt$-modules. Then by definition
$[T_k,T_1]_D=[T_k,T_1]_{k\Qalt}$ and moreover
$[T_k,T_1]_{k\Qalt}=[\overline{T_k},\overline{T_1}]_{k\Qalt\setminus\{l\}}$
 because
$l$ is not in the support of $T_1$ nor $T_k$. But this is zero since
$[\overline{T_k},\overline{T_1}]_{\overline{\C}}=0$.
\item $T_k$ is a $k\Qalt$-module  and
 $T_1=P_h[1]$  for some $h\ne l$. Then by definition
 $[T_k,T_1]_D=[T_k,P_h]^1_{k\Qalt}$ and
$[T_k,P_h]^1_{k\Qalt}
=[\overline{T_k},\overline{P_h}]^1_{k\Qalt\setminus\{l\}}$
 because
$l$ is not in the support of $T_1$ nor $T_k$. But this is zero since
$[\overline{T_k},\overline{T_1}]_{\overline{\C}}=0$.
\item  $T_k$ is the first shift of a projective
 and $T_1$ is a $k\Qalt$-module then by definition $[T_k,T_1]_D=0$.
\item $T_1$ and $T_k$ are both shifts of projectives, say
 $T_1=P_h[1]$ and $T_k=P_{h'}[1]$, then $[T_k,T_1]_\C=1$ implies that
 $h\to h'$ is an arrow in $\Qalt$, since $\Qalt$ is alternating.
 Moreover $T_k\to T_1$ is an arrow in the Auslander-Reiten quiver of
 $\C$ and hence $T_1\to T_k$ is an arrow in $\QT$. By hypothesis, we
 also have an arrow $T_1\ot T_k$ in $\QT$. This is impossible, since
 $\QT $  is a cluster diagram by Theorem \ref{quivers}.
\end{enumerate}
Thus $[T_k,T_1]_D=0$.
Moreover, it follows from the fact that $T_1,T_k$ are indecomposable
modules or first shifts of projectives,
that $[T_k,\tau T_1[-1]]_D=0$.
Hence
\begin{eqnarray}
1=[T_k,T_1]_\C&=&
[T_k,\tau T_1[-1]]_D+[T_k, T_1]_D+[T_k,\tau^{-1}T_1[1]]_D\nonumber\\
&=&[T_k,\tau^{-1}T_1[1]]_D \label{eq20}\\
&=& [T_k,\tau^{-1}T_1]^1_D,\label{eq21}\end{eqnarray}
and also
$$1=[\tau T_k,T_1]^1_D.$$
By calculations similar to those preceding (\ref{eq20},\ref{eq21}), one
can show that
$$0=[\overline{T_k},\overline{T_1}]_{\overline{\C}}=
[\overline{\tau T_k},\overline{T_1}]^1_D$$
Since $T_1$
is $l$-free, the restriction does not change $T_1$. Thus the restriction
must change $\tau T_k$ and hence $\tau T_k$ is not $l$-free.
Here we use the fact that if two indecomposable objects $M,N$ are
$l$-free then $[\overline{M},\overline{N}]^1_D=[{M},N]^1_D$.
A similar argument shows that $\tau^{-1} T_1$ is not $l$-free.
That is
$1=[P_l[1],\tau T_k]^1_\C=[T_k,P_l[1]]_\C$ and
$1=[\tau^{-1}T_1,P_l[1]]^1_\C=[P_l[1],T_1]_\C$.
But then there
 are two paths $q_1,q_2$ in $\QT$, $q_1$ going from $T_1$ to $P_l[1]$
 and $q_2$ from $P_l[1]$ to $T_k$, and $q_1,q_2$
 are both non-zero  in $(\QT,I_T)$.
The composition $q=q_1q_2$ is a path from $T_1$ to $T_k$.
This path $q$ is not a shortest path because of the hypothesis $m=1$ and
by Proposition \ref{shortest}(3), we have $q=0$ in $(\QT,I_T)$.
Consider the lifts $q^D_2\in\Hom_D(T_k,P_l[1])$ and $q_1^D\in
\Hom(P_l[1],\tau^{-1}T_1[1])$. Recall our convention that $T_k,T_1$ are
$k\Qalt$-modules or first shifts of projectives.
Now $q$ being zero in $(\QT,I_T)$ means that the composition
$q^D_1\circ q^D_2=0$, whence
$[T_k,\tau^{-1}T_1[1]]_D=0$, contradiction to (\ref{eq20}).

We have shown that the path $p:T_1\to T_2 \to \ldots\to T_k$ passes through
all vertices of $\QT$, that is $k=n$.
Since  $p$ is a shortest path, the underlying graph of the quiver $\QT$ is
a cycle. Thus $\QC=\QT=T_1\to T_2 \to\ldots\to T_n \to T_1$ with $n>3$
 and by a result
of \cite{cluster2} this implies that the cluster algebra (and hence
$\Qalt$) is of type $D_n$.
Let us label the  vertices of $\Qalt$ as follows.
\begin{diagram}[size=2em]
&&&&&&&&&&n-1 \\
&&&&&&&&&\ruLine \\
1&\rLine&2&\rLine&3&\rDots&(n-3)&\rLine& (n-2)&\\
&&&&&&&&&\rdLine \\
&&&&&&&&&&n \\
\end{diagram}
Note that if one removes any vertex $T_i$ of the quiver $\mathcal{Q}_T$ then
the induced subquiver $\mathcal{Q}_T-\{T_i\}$ is a Dynkin quiver of type
$A_{n-1}$. Then the corresponding cluster category
$\C_{\mathcal{Q}_T-\{T_i\}}$ is of type $A_{n-1}$ too. Therefore
the position of $T_i$ in the Auslander-Reiten quiver of $\C$ must be at
level $n$ or $n-1$; that is, $T_i=\tau^{-k}(P_l)$ with $l\in\{n-1,n\}$
and $k\ge 0$.
Suppose without loss of generality that $T_1=P_l[1]$.
Now since $T$ is a tilting object and since there are
arrows $T_i\to T_{i+1}$ in $\mathcal{Q}_T$, we have for all $2\le i\le
n-1$
$$T_i=\left\{\begin{array}{ll}
\tau^{i-1}P_l[1] &\textup{if $i$ is odd}\\
\tau^{i-1}P_{l'}[1] &\textup{if $i$ is even}
\end{array}\right.
\quad\textup{and}\quad
  T_{n}\ =\ P_{l'} $$
where $l'$ is such that $\{l,l'\}=\{n-1,n\}$. Thus $[T_n,T_1]_\C=0$ and
consequently $p$ is zero in $(\QT,\IT)$.

Suppose now that $m=2$. By Lemma \ref{hom} we have $[T_j,T_i]_\C\le 1$ and
therefore either $p_1=p_2$ in $(\QT,\IT) $ (and in this case we are
done) or one of $p_1,p_2$, say $p_1$
is zero in $(\QT,\IT)$ and $p_2$ is not zero. Thus $[T_j,T_i]_\C=1$ and
then
$p_1$ being zero means that there is a vertex $T_h$ on the path $p_1$ such
that
\begin{equation}\label{eq13}[T_j,T_h]_\C=0.
  \end{equation}
We may suppose
that $T_h=P_l[1]$ for some $l$.
As we did before in the case $m=1$, we remove that
vertex $T_h$ so that we get a quiver $(\qt,I_{\overline{T}})$ of rank $n-1$. In this quiver,
the induced path $\overline{p_2}$ is zero by case $m=1$. We have seen in
the case $m=1$ that if $p_2$ is non-zero in $(\QT,\IT)$ then
$[T_j,P_l[1]]_\C=1$, contradiction to (\ref{eq13}).
\end{proof}

\end{section}


\begin{section}{Applications}\label{appl}
In this section, we use the results in section \ref{qar} to answer
conjectures of \cite{BMRRT} and \cite{CCS}.
We keep the setup of the previous section.
Denote
by $\nu $ the number of positive roots of the root system corresponding to
the type of the cluster algebra.
By a result of \cite{BMR}, the number of indecomposable
$\End_\C(T)^{\textup{op}}$-modules is equal to $\nu$.
On the other hand, in \cite{CCS} it has been shown for type $A$
(and conjectured for types $D$ and $E$) that the number of
indecomposable $k\QC/\langle I_C\rangle$-modules is also equal to $\nu$.
Using this, we can prove in type $A$ the following theorem, which
has been conjectured in \cite{BMRRT} for types $A,D,E$.
\begin{theorem}\label{A}
 Let $C$ be any cluster of a cluster
algebra of type $A$ and let $(\QC,I_C)$ be its quiver with relations.
Let $T$ be a corresponding tilting object in the cluster category.
Then the  cluster tilted algebra
$\End_\C(T)^{\textup{op}}$ is isomorphic to the algebra $k \QC/\langle I_C\rangle$.
\end{theorem}
Using Theorem \ref{quivers}, Proposition \ref{relations} and the
considerations above, the result follows from
\begin{lemma}
Let  $A$ be an algebra and let $I$ be an ideal of $A$. Suppose that the
category mod$A$ of finitely generated $A$-modules has a finite number of
indecomposable modules. Suppose also that the category mod$A/I$ has the
same number of indecomposable modules. Then, $I$ is zero.
\end{lemma}
\begin{proof}
 We denote by $\overline{\textup{mod}} A$ the category of isoclasses of
$A$ modules. Let $j$ be the natural map from mod$A/I$ to mod $A$.  It is
clear that $j$ gives a quotient map from $\overline{\textup{mod}} A/I$
to $\overline{\textup{mod}} A$. We still denote this map by $j$. The
image of $j$ is the subcategory of isoclasses of $A$-modules on which $I$
vanishes. Moreover, $j$ commutes with direct sums, hence, it sends
indecomposable modules on indecomposable ones. It is easily seen that $j$
is injective, so $j$ embeds the set of isoclasses of indecomposable $A/I$
modules in the set of isoclasses of indecomposable $A$ modules. By the
hypothesis of the lemma, this restriction of $j$ is bijective. Hence, by
the Krull-Schmidt theorem, $j$ is bijective. This implies that $I$
vanishes on all finitely generated $A$-modules. Considering $A$ as an
$A$-module then gives $I=0$.
\end{proof}

Next, we describe the exchange relations of the cluster algebra
in terms of the cluster category.
Let $M$ be an indecomposable summand of $T$ and
$T=\overline{T}\oplus M$ with $\overline{ T}=\oplus_{i=1}^{n-1}T_i$.
Suppose that $M'$ is another indecomposable
object of the cluster category such that $M$ and $M'$ form an exchange
pair, that is $T'=\overline{T}\oplus M'$ is another tilting object. By
\cite{BMRRT}, there are two triangles in $\C$
$$\begin{array}{cccccccccc}
M & \to & \oplus_{i\in I'} T_i & \to & M' &\to & M[1]\\
M'&\to &\oplus_{i\in I} T_i&\to&M&\to&M'[1]
\end{array}
$$
Let $z,z',x_i$ be the cluster variables corresponding to $M,M',T_i$
respectively and let $C,C'$ the clusters corresponding to $T,T'$.
Let $\mathcal{B}=(b_{xy})_{x,y\in C} $ be the sign-skew-symmetric matrix
associated to the cluster $C$, see \cite{cluster2}.
Then
$C'=C\setminus\{z\}\cup\{z'\}$, $C\cap C'=\{x_1,\ldots,x_{n-1}\}$
and $z,z'$ satisfy the so called
exchange relation:
\begin{equation}\label{eq12}
zz'=\prod_{x\in C : b_{xz}>0}x^{b_{xz}}+\prod_{x\in C :
b_{xz}<0}x^{-b_{xz}}.
\end{equation}
The following theorem has been conjectured in \cite{BMRRT}.
\begin{theorem}\label{exchange} For any cluster algebra of type $A,D,E$,
 in the situation above the exchange relation can be written as
 $$zz'= \prod_{i\in I} x_i+\prod_{i \in I'}x_i.$$
\end{theorem}
\begin{proof}
  By definition of the cluster diagram $\QC$, there is a vertex for each
  cluster variable $x$ in $C$ and there is an arrow $x\to y$ precisely
  if $b_{xy}>0$. Since the cluster algebra  is  of type $A,D$ or $E$,
  we have $b_{xy}\in \{-1,0,1\}$. Thus (\ref{eq12}) becomes
  $$zz'=\prod_{x\to z \textup{ in } \QC} x +\prod_{x\ot z \textup{ in
  }\QC} x. $$
  Now the result follows from Theorem \ref{quivers} and Lemma \ref{BB'}.
\end{proof}

Finally, we generalize a result of \cite{CCS} on denominators of Laurent
polynomials. This theorem has been conjectured in \cite{CCS} in a slightly
different form using the quiver with relations $(\QC,I_C)$ instead of
the cluster category.
\begin{theorem}\label{denominator}
  Let $C=\{x_1,\ldots,x_n\}$ be any cluster of a cluster algebra of type
  $A,D$ or $E$ and let $T=\oplus_{i=1}^n T_i$ the corresponding tilting
  object in the cluster category $\C$. Then there is a bijection
  $$\begin{array}{ccc}
  \{\textup{indecomposable objects of $\C$}\}&\to &\{\textup{cluster
  variables}\}\\
  M&\mapsto&x_M\\
  \end{array}$$
  such that
  \begin{equation}\label{eq13a}
  x_M = \frac{P(x_1,\ldots,x_n)}{\prod_{i=1}^n
  x_i^{[T_i,M]_\C^1}}\end{equation}
  where $P$ is a polynomial prime to $x_i$ for all $i$ and
  $[T_i,M]_\C^1=\dim\Ext_\C(T_i,M).$
  \end{theorem}
  \begin{remark}
    It has been shown in \cite{BMR} that $\Hom_\C(\tau^{-1}T,\ )$ induces
    an equivalence of categories $\C/\textup{add} T \to
    \textup{mod}\End_\C(T)^{\textup{op}}$.
    Under this equivalence, the object $M$ of  $\C$ gets mapped to
    the indecomposable $\End_\C(T)^{\textup{op}}$-module with dimension
    vector $(d_1,\ldots,d_n)$, $d_i=\dim\Ext(T_i,M)$. Thus the exponent
    of $x_i$ in the denominator is the multiplicity of the simple
    $\End_\C{T}^{\textup{op}}$-module $S_i$ in the image of $M$.
  \end{remark}
\begin{proof}
  The existence of the bijection between the two sets is proved in
  \cite{BMRRT}. The fact that $x_M$ can be written in terms of the
  $x_1,\ldots,x_n$ as a Laurent polynomial is the Laurent phenomenon
  proved in \cite{cluster1}. We have to show that the exponents in the
  denominator are as stated.
  It has been shown in \cite{cluster2} that cluster variables are in
  bijection with almost positive roots. Let $\za_M,\za_i$ be the almost
  positive root corresponding to $x_M,x_i$ respectively. By a result of
  \cite[Prop. 6.5]{CCS}, the exponent of $x_i$ in the denominator of (\ref{eq13a})
  is equal to the compatibility degree $(\za_i\mid\mid\za_M)$ of the
  almost positive roots. Finally, the identity $(\za_i\mid\mid \za_M)=
  \dim\Ext_\C(T_i,M)$ has been shown
  in \cite{BMRRT}.
  \end{proof}

\begin{example}
  We give an example of type $D_5$. Let $C$ be a cluster having the
  following diagram $\QC$
  \begin{diagram}[size=2em]
    5&\rTo     &&&1&\rTo &&& 2\\
     &\luTo&&\ldTo&&\luTo &&\ldTo\\
     &&     4      &\rTo&&&3
  \end{diagram}
  Performing a mutation at vertex $1$ followed by a mutation at vertex
  $2$ shows that $\QC$ is mutation equivalent to a quiver with
  underlying graph the Dynkin diagram $D_5$. A corresponding tilting
  object $T=\oplus_{i=1}^n T_i $ is illustrated in the Auslander-Reiten
  quiver of $\C$ in figure \ref{example},
  where $T_i$ corresponds to the vertex $i$ of $\QC$. Let $M$ be the
  indecomposable  object shown in the same figure.
  \begin{figure}
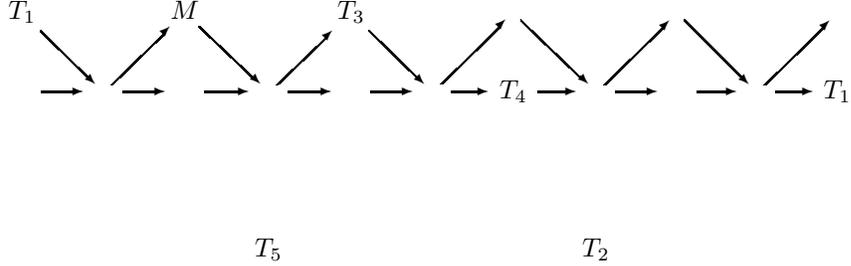

    \begin{diagram}[size=1.5em]
      &&T_1&&&&M&&&&T_3&&&& &&&& &&&& \\
      &&&\rdTo&&\ruTo&&\rdTo&& \ruTo&&\rdTo&&\ruTo&&\rdTo&&
\ruTo&&\rdTo&&\ruTo&&&&        \\
      &&\ \ &\rTo&\ \ &\rTo&\ \ &\rTo&\ \ &\rTo&\ \ &\rTo&\ \
      &\rTo&T_4&\rTo&\   \
      &\rTo&\ \ &   \rTo&\ \ &\rTo&T_1&&\\
     &&&\ruTo&&\rdTo&& \ruTo&&\rdTo&&\ruTo&&\rdTo&&
\ruTo&&\rdTo&&\ruTo&&\rdTo&&        \\
     \\
      &&&\rdTo&&\ruTo&&\rdTo&& \ruTo&&\rdTo&&\ruTo&&\rdTo&&
\ruTo&&\rdTo&&\ruTo&&        \\
       &&&& &&&& T_5&&&&&&&&T_2&&&&\\
    \end{diagram}
    \caption{Auslander-Reiten quiver with tilting object}\label{example}
  \end{figure}
  By Theorem
  \ref{denominator} we have
  $$x_M=\frac{P(x_1,\ldots,x_n)}{x_1x_2x_3x_4}.$$
  Note that the shape of the Auslander-Reiten quiver of the cluster tilted algebra is
  obtained from figure \ref{example} by deleting the vertices
  $T_1,\ldots,T_5$. At the position of $M$ we find the indecomposable
  $\End_\C(T)^{\textup{op}}$-module $\Hom_\C(\tau^{-1} T,M)$, by
  \cite{BMR}. It is the first projective and the third injective
  indecomposable  of $\End_\C(T)^{\textup{op}}$.
  Its dimension
    vector $\mathbf{d}=(d_1,d_2,d_3,d_4,d_5)$ is given by
    $d_i=[\tau^{-1} T_i,M]_\C=[M,T_i]_\C^1$, thus $\mathbf{d} =
    (1,1,1,1,0)$.
\end{example}

\end{section}

\appendix
\begin{section}{Appendix}

In this section, we give some general results about cluster quivers of
simply-laced finite type.

\begin{theorem}
  For each cluster quiver $Q$ of simply-laced finite type, there
  exists an embedding of $Q$ in the plane such that all vertices of
  $Q$ belong to the closure of the unique unbounded connected
  component of the complement.
\end{theorem}

Let us call this a \textit{nice embedding}. It is clear that each
bounded connected component is a topological cell.

\begin{corollary}
  In a nice embedding, the boundary of each bounded cell is an
  oriented cycle in the quiver. There are no other edges between its
  vertices.
\end{corollary}
\begin{proof}
  It is enough to prove that there are no other edges between the
  boundary vertices of the fixed cell, as this is known to imply the
  orientation property \cite[Prop. 9.7]{cluster2}. As there are no
  edges inside the cell, other edges must be outside. Would they
  exist, they would contradict the property of the embedding that all
  vertices border the unbounded component.
\end{proof}

\begin{corollary}
  For each arrow $i \to j$ in $Q$, there are at most two shortest
  paths from $j$ to $i$. The shortest paths from $j$ to $i$ are given
  by the boundaries of the bounded cells which are adjacent to the
  arrow $i \to j$ in any nice embedding.
\end{corollary}
\begin{proof}
  Either side of the arrow $i \to j$ is either a bounded cell or an
  unbounded component. If one of those sides is a bounded cell, it
  gives a shortest path from $j$ to $i$.

  Conversely, pick a shortest path from $j$ to $i$. By definition,
  there is no other edge between its set of vertices. As there can be
  no other vertex inside the loop drawn by the shortest path (this
  would contradict the nice embedding property), one deduces that this
  loop bounds a cell, which is of course adjacent to the arrow $i \to
  j$.
\end{proof}

Let us now prove the theorem.
This could be done by inspection of all possible cluster quivers
of simply-laced finite type but we use another proof.

\begin{proof}
  The strategy of proof is the following one. First the Theorem is
  clearly true for oriented Dynkin diagrams, which are trees. Any
  plane embedding of a tree is nice. We are going to prove that the
  statement of the Theorem is stable by mutation of quivers. This is
  clearly true if the mutation does not change the shape of the
  quiver.

  To go further, it is necessary to have a precise description of what
  can happen during the mutation process. We need to describe all
  possible configurations around a vertex of a cluster quiver of
  simply-laced finite type.

  Recall that the \textit{link} of a vertex $v$ in a graph is the
  graph induced on the set of vertices which are adjacent to $v$.

  Let us start with some simple remarks on the link of a vertex in a
  cluster quiver of simply-laced finite type.

  First, the link of any vertex has at most 3 connected components.
  This follows from the fact that no orientation of the affine $D_4$
  diagram is of finite type.
  
  Next, as each triangle with vertex $v$ must be oriented, each vertex of
  the link is either a sink or a source in the link.

  We claim that each connected component of the link is a tree.
  Indeed, there can not be any odd cycle, because sources and sinks
  must alternate. It is also easy to check that the existence of a
  4-cycle or a 6-cycle would imply that the quiver is not of finite
  type. Any even cycle of length at least 8 would imply that the
  quiver contains an affine $D_4$ quiver, which is not of finite type.

  For similar reasons, each connected component of the link is a
  linear tree. It is enough to prove that the existence of a fork
  would contradict the assumption that the quiver is of finite type.
  This is readily checked by a sequence of mutations.

  Hence we know that each connected component of the link is an
  alternating linear tree. So we can describe each connected component
  by its cardinality, up to reversal, and a link can be described by
  the set of cardinalities of its connected components.

  Here is a list of all possible links:
  \begin{itemize}
  \item One connected component: (1);(2);(3);(4);(5);(6),
  \item Two connected components: (1,1);(1,2);(2,2);(1,3);(1,4);(2,3);(2,4),
  \item Three connected components: (1,1,1);(1,1,2);(1,2,2).
  \end{itemize}
  
  Indeed the links $(7)$;$(3,3)$;$(1,5)$ and $(1,1,3)$ are not
  possible as they would give that an affine $D_4$ quiver is of finite
  type. Similarly the link $(2,2,2)$ contains an affine $E_6$ quiver.

  For each link, there is only one possible orientation up to global
  change of orientation, except for $(1,1,2)$ and $(1,3)$ where there
  are two really different orientations.

  Conversely each of those links are realized in a cluster quiver of
  finite type.
 
  Here comes now the list of all possible shape-changing mutations:
  \begin{align}
    (2) &\lra (1,1),\\
    (3) &\lra (1,1,1),\\
    (5) &\lra (1,2,2),\\
    (6) &\lra (2,4),\\
    (4) &\lra (1,1,2),\\
    (1,3) &\lra (1,1,2),\\
    (1,2) &\lra (1,2),\\
    (2,2) &\lra (2,2),\\
    (1,3) &\lra (1,3),\\
    (2,3) &\lra (2,3).
  \end{align}

  To prove the Theorem, one now has to check in both directions for
  each of these cases that the existence of a nice embedding before
  mutation permits to build a nice embedding after mutation.

  The principle is the same for all cases. Pick one of these links and
  assume it is part of a nice embedding. The finiteness assumption and
  the nice embedding hypothesis together allow to give restrictions on
  the local picture of the embedding near the fixed vertex $v$. Then
  using these restrictions, one concludes that the quiver after
  mutation still has a nice embedding.

  In general, one knows that at least one of the components near $v$
  not enclosed by the link of $v$ has to be unbounded.

  This is enough to solve the cases $(1,1)\lra(2)$, $(1,1,1)\lra(3)$ and
  $(1,2)\lra(1,2)$.

  CLAIM : let $i \to j$ be an arrow in the link of $v$. Assume that
  there is a bounded cell containing this arrow but not $v$. Then the
  only edges between the vertices of this cycle (but $i$ and $j$) and
  a vertex of the link of $v$ are the edges from $i$ or $j$ to their
  neighbour in the cycle.

  Indeed, the existence of such a vertex and edge would contradict the
  nice embedding property.

  CLAIM : replacing each cell containing an arrow of the link of $v$
  and not containing $v$ by a triangle gives a quiver of finite type.

  Indeed one can show that all these cycles can only meet or be
  related by an edge inside the link of $v$. Hence one can use
  mutation to shorten the cycles independently until they become
  triangles.

  In some cases, it is necessary to show that at least one of some
  arrows in the link has an unbounded side. By the claim above, this
  is done by checking that adding triangles on all these arrows can
  not give a quiver of finite type. Then one can repeat this argument
  to get more information on the local configuration.

  Combined with the argument on the unbounded cell near $v$, this is
  enough to solve all the remaining cases.

  Figure \ref{paires} displays all possible plane configurations
  around a vertex $v$ in a nice embedding. The unbounded component is
  shaded and the vertex $v$ is marked. Each configuration is paired
  with the configuration obtained after mutation.

  \begin{figure}
    \begin{center}
      \scalebox{0.5}{\includegraphics{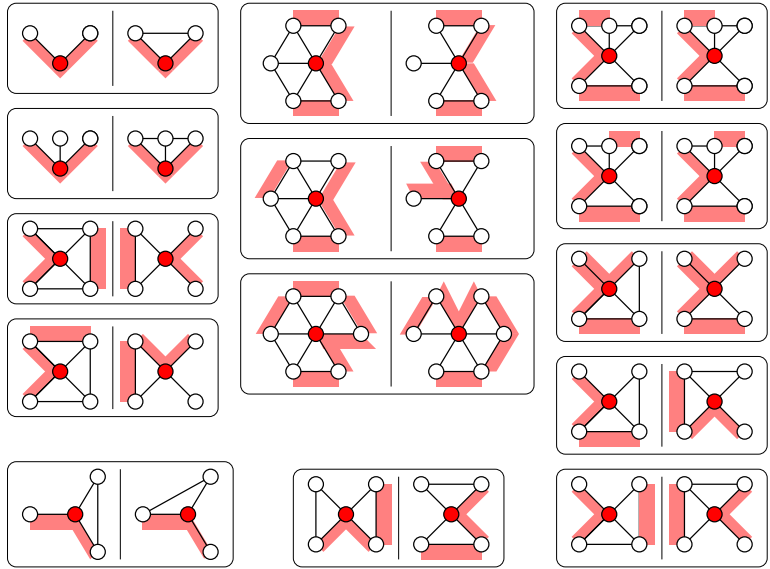}}
      \caption{List of all local mutation configurations}
      \label{paires}
    \end{center}
  \end{figure}

\end{proof}
\end{section}

\newcommand{\etalchar}[1]{$^{#1}$}


\end{document}